\documentclass{article}
\usepackage[english]{babel}
\usepackage{amsmath}
\usepackage{amsfonts}
\usepackage[latin1]{inputenc}
\usepackage[T1]{fontenc}
\usepackage{amssymb}
\usepackage{color}
\usepackage{cite}
\usepackage{amsthm}
\usepackage{latexsym}
\usepackage{graphicx}
\usepackage[]{geometry}
\usepackage[ruled,algosection,linesnumbered]{algorithm2e}
\newcommand{\Mu}{\boldsymbol{\mu}} 
\usepackage[lofdepth,lotdepth]{subfig}

\newtheorem{theorem}{Theorem}[section]

\newtheorem{proposition}[theorem]{Proposition}

\newtheorem{remark}[theorem]{Remark}

\title{Static condensation optimal port/interface reduction and error estimation for structural health monitoring}
\author{Kathrin Smetana\footnote{Department of Applied Mathematics, University of Twente, P.O. Box 217, 7500 AE Enschede, The Netherlands, \texttt{k.smetana@utwente.nl}; also: Department of Mechanical Engineering, Massachusetts Institute of Technology, Cambridge, MA
02139, United States.}}

\begin{document}
\maketitle
\begin{abstract}
Having the application in structural health monitoring in mind, we propose reduced port spaces that exhibit an exponential convergence for static condensation procedures on structures with changing geometries for instance induced by newly detected defects. Those reduced port spaces generalize the port spaces introduced in [K. Smetana and A.T. Patera, SIAM J. Sci. Comput., 2016] to geometry changes and are optimal in the sense that they minimize the approximation error among all port spaces of the same dimension. Moreover, we show numerically that we can reuse port spaces that are constructed on a certain geometry also for the static condensation approximation on a significantly different geometry, making the optimal port spaces well suited for use in structural health monitoring.
\end{abstract}


\section{Introduction}
\label{sec:1}

Manual or automated inspection of large structures such as offshore platforms is carried out on a regular basis; the effects of any detected defects must be assessed rapidly in order to avoid further damage or even catastrophic failure. This can be facilitated by relying on numerical simulations. One step towards a fast numerical simulation response for such large structures is to exploit their natural decomposition into components and apply static condensation to obtain a (Schur complement) system of the size of the degrees of freedom (DOFs) on all interfaces or ports in the system.  However, as the size of this Schur complement system may still be very large it is vital to reduce the number of DOFs on the interfaces or ports and thus consider reduced interface or port spaces. 

In the popular component mode synthesis (CMS) approach \cite{Hur65,BamCra68,Bou92,HetLeh10,JaBeLa11} this reduced space is spanned via certain eigenmodes. In \cite{HuKnPa13} generalized Legendre polynomials are used. In \cite{FeHoEb16} deformation patterns from an analysis of the assembled structure are employed. Moreover, local reduced models are generated from parametrized Lagrange or Fourier modes and coupled via FE basis functions in \cite{IaQuRo16}. Finally, empirical modes generated from local solutions of the PDE are suggested in \cite{EftPat13a,MarRozHaa14,Betal17}.

Recently, port spaces that are optimal in the sense of Kolmogorov and thus minimize the approximation error among all port spaces of the same dimension have been introduced in \cite{SmePat16}. The approach in \cite{SmePat16} generalizes the idea of separation of variables by connecting two components at the port for which we wish to construct the port space and consider the space of all local solutions of the partial differential equation (PDE) with arbitrary Dirichlet boundary conditions on the ports that lie on the boundary of the two-component system. From separation of variables we anticipate an exponential decay (of the higher modes) of the Dirichlet boundary conditions to the interior of the system. To quantify which information of the Dirichlet boundary conditions reaches the shared port of the system, a (compact) transfer operator that acts on the space of local solutions of the PDE is introduced. Solving the transfer eigenproblem for the composition of the transfer operator and its adjoint yields the optimal space. For related work in the context of the generalized finite element method we refer to \cite{BabLip11,BaHuLi14}. 

In \cite{BuhSme18} it has been shown that by employing methods from randomized numerical linear algebra an extremely accurate approximation of those optimal port spaces can be computed in close to optimal computational complexity. To account for variations in a material or geometric parameter in \cite{SmePat16} a \emph{parameter-independent} port space is generated from the optimal parameter-dependent port spaces via a spectral greedy algorithm. It is further numerically demonstrated in \cite{SmePat16} that the optimal port spaces often outperform other approaches such as Legendre polynomials \cite{HuKnPa13} or empirical modes \cite{EftPat13a}; also an exponential convergence can be observed. Finally, those optimal port modes have been used in structural integrity management of offshore structures in \cite{Ketal18} and optimal local approximation have been exploited in the context of data assimilation in \cite{TadPat18}.

In this article we want to investigate the applicability of optimal port spaces for structural health monitoring and more specifically extend the concept of \cite{SmePat16} to geometry changes. First, we show how to construct one port space for several different geometries such as a beam and a beam with a crack or hole via a spectral greedy algorithm. Moreover, if during an inspection a defect is detected, unfortunately, often the precise geometry of the newly detected defect is not amongst the component geometries the reduced model has been trained for. Therefore, we demonstrate numerically that for realistic error tolerances the optimal port spaces constructed on one geometry can often be reused on another. In order to assess whether the resulting reduced model is accurate enough we suggest to employ the error estimator for port reduction introduced in \cite{Sme14} as this error estimator is both an upper and lower bound of the error and based on local error indicators associated with the ports; also the latter are a lower bound of the local error on the component pair that shares the respective port. Error estimation for port or interface reduction has also been considered in \cite{EftPat13a,Betal17,JaBeLa11}. Finally, we note that also in \cite{Betal17} local reduced order models for geometry changes are suggested. However, the authors of \cite{Betal17} neither reuse existing reduced models nor build one reduced model for different geometries. 

The remainder of this paper is organized as follows. In section \ref{sec:2} we introduce the problem setting and recall the algebraic (port reduced) static condensation procedure. Subsequently, we recall the optimal port spaces introduced in \cite{SmePat16} in section \ref{sec:3}. In section \ref{sec:4} we propose quasi-optimal port spaces for parametrized problems including geometry changes such as from a beam to a beam with a crack. Subsequently, we discuss in section \ref{sec:5} how to deal with systems with many components and recall in section \ref{sec:6} the {\it a posteriori} error estimator from \cite{Sme14}. Finally, we present numerical experiments in section \ref{sec:7} and draw some conclusions in section \ref{sec:8}.

\section{Preliminaries}
\label{sec:2}

\subsection{Problem setting}

Let $\Omega_{gl} \subset \mathbb{R}^{d}$, $d=2,3$, be a large, bounded domain with Lipschitz boundary and assume that $\partial \Omega_{gl}= \Sigma_{D} \cup \Sigma_{N}$, where $\Sigma_{D}$ denotes the Dirichlet and $\Sigma_{N}$ the Neumann boundary, respectively. We consider a linear, elliptic PDE on $\Omega_{gl}$ with solution $u_{gl}$, where $u_{gl}$ equals $g_{D}$ on $\Sigma_{D}$ and satisfies homogeneous Neumann boundary conditions on $\Sigma_{N}$ noting that the extension to non-homogenoues Neumann boundary conditions is straightforward. 

\begin{figure}[t]
\begin{center}
\includegraphics[scale = 0.12]{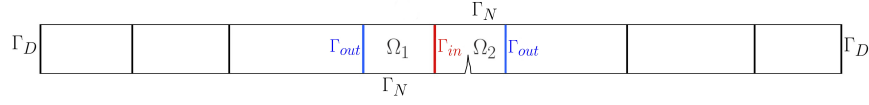}
\caption{Illustration of $\Omega_{gl}$, $\Omega$, defined as $\bar{\Omega} = \bar{\Omega}_{1} \cup \bar{\Omega}_{2}$, and the ports $\Gamma_{in}$ and $\Gamma_{out}$.}\label{fig:omega}
\end{center}
\end{figure}

To compute an approximation of $u_{gl}$ we decompose the large domain $\Omega_{gl}$ into (many) non-overlapping subdomains. To simplify the presentation we consider henceforth two subdomains $\Omega_{1}, \Omega_{2} \subset \Omega_{gl}$ and their union $\Omega$ with $\bar{\Omega} = \bar{\Omega}_{1} \cup \bar{\Omega}_{2}$ as illustrated in Fig.~\ref{fig:omega}; the approximation of the whole system associated with $\Omega_{gl}$ will be discussed in Sec.~\ref{sec:5}. Moreover, we introduce the shared interface $\Gamma_{in}:=\bar{\Omega}_{1} \cap \bar{\Omega}_{2}$ and $\Gamma_{out}:=\partial\Omega \setminus \partial \Omega_{gl}$. 

We consider the following problem on $\Omega$: For given $f \in L^{2}(\Omega)$ find $u$ such that
\begin{equation}\label{eq:PDE}
\mathcal{A}u = f \enspace \text{in } \Omega, \quad \text{and} \quad u = u_{gl} \enspace \text{on } \Gamma_{out},
\end{equation}
where $\mathcal{A}$ is a linear, elliptic, and continuous differential operator. We may then introduce a conforming Finite Element (FE) discretization and a FE approximation $\mathsf{u}$ whose FE coefficients $\underline{\mathsf{u}} \in \mathbb{R}^{N}$ solve the following  linear system of equations
\begin{equation}\label{eq:FEM}
\underline{A} \underline{\mathsf{u}} = \underline{\mathsf{f}}.
\end{equation}
Here, $\underline{A} \in \mathbb{R}^{N\times N}$ discretizes the (weak form of the) differential operator $\mathcal{A}$ and $\underline{\mathsf{f}} \in \mathbb{R}^{N}$ accounts for the discretization both of $f$ and enforcing the non-homogeneous Dirichlet boundary conditions $u_{gl}|_{\Gamma_{out}}$; we assume that in the rows associated with the Dirichlet DOFs the non-diagonal entries are zero and the diagonal entries equal one.

\subsection{Static condensation}\label{subsec:static_cond}
To obtain a linear system of equations of the size of the number of DOFs $N_{in}$ on the interface $\Gamma_{in}$ we perform static condensation. To that end, we first sort the DOFs in DOFs associated with $\Omega_{1}$, $\Omega_{2}$, and $\Gamma_{in}$ to rewrite \eqref{eq:FEM} as follows:
\begin{equation}
\begin{bmatrix}
   \underline{A}_{\Gamma_{in}} &  \underline{A}_{\Gamma_{in},\Omega_1}^{T} & \underline{A}_{\Gamma_{in},\Omega_2}^{T}  \\
    \underline{A}_{\Gamma_{in},\Omega_1} & \underline{A}_{\Omega_{1}} & 0 \\
     \underline{A}_{\Gamma_{in},\Omega_2} & 0 & \underline{A}_{\Omega_{2}} 
\end{bmatrix}
\begin{bmatrix}
\underline{\mathsf{u}}_{\Gamma_{in}}\\
\underline{\mathsf{u}}_{\Omega_{1}} \\
\underline{\mathsf{u}}_{\Omega_{2}}
\end{bmatrix}
= 
\begin{bmatrix}
\underline{\mathsf{f}}_{\Gamma_{in}}\\
\underline{\mathsf{f}}_{\Omega_{1}} \\
\underline{\mathsf{f}}_{\Omega_{2}}
\end{bmatrix}.
\end{equation}
We may then apply static condensation to remove the DOFs corresponding to the interior of $\Omega
_{1}$ and $\Omega_{2}$: We define the Schur complement matrix and the Schur complement right-hand as
\begin{align}\label{eq:SC_matrix}
\underline{A}_{SC} &= \underline{A}_{\Gamma_{in}} - \underline{A}_{\Gamma_{in},\Omega_1}^{T} \underline{A}_{\Omega_{1}}^{-1}  \underline{A}_{\Gamma_{in},\Omega_1} - \underline{A}_{\Gamma_{in},\Omega_2}^{T}\underline{A}_{\Omega_{2}}^{-1} \underline{A}_{\Gamma_{in},\Omega_2} \quad \in \mathbb{R}^{N_{in} \times N_{in}}\\
\underline{\mathsf{f}}_{SC} &= \underline{\mathsf{f}}_{\Gamma_{in}} - \underline{A}_{\Gamma_{in},\Omega_1}^{T} \underline{A}_{\Omega_{1}}^{-1}\underline{\mathsf{f}}_{\Omega_{1}} -  \underline{A}_{\Gamma_{in},\Omega_2}^{T}\underline{A}_{\Omega_{2}}^{-1}\underline{\mathsf{f}}_{\Omega_{2}} \quad \in \mathbb{R}^{N_{in}}
\end{align}
such that the vector of interface coefficients solves the \emph{Schur complement system} 
\begin{equation}\label{eq:SC}
\underline{A}_{SC} \underline{\mathsf{u}}_{\Gamma_{in}} = \underline{\mathsf{f}}_{SC} \qquad \text{of size } N_{in}\times N_{in}.
\end{equation}
We note that computing $\underline{A}_{\Omega_{i}}^{-1}  \underline{A}_{\Gamma_{in},\Omega_i}$ corresponds to solving the PDE on $\Omega_{i}$, $i=1,2$, $N_{in}$ times with homogeneous Dirichlet boundary conditions on $\Gamma_{in}$ and right-hand sides that occur from lifting the respective $N_{in}$ FE basis functions on $\Gamma_{in}$. 

\subsection{Port or interface reduction}\label{subsect:port red}

As indicated in the introduction, unfortunately, for many real-world applications the size of \eqref{eq:SC} is still too large, such that a further reduction in size is desirable. We assume that we have a reduced basis $\underline{\phi}_{1}, \hdots, \underline{\phi}_{n} \in \mathbb{R}^{N_{in}}$, $n \ll N_{in}$, at our disposal, which we store in the columns of a matrix $\underline{\Phi}_{n} \in \mathbb{R}^{N_{in} \times n}$. We may then introduce a \emph{port reduced static condensation approximation} \cite{EftPat13a} $\mathsf{u}^{n}$ with FE coefficients $\underline{\mathsf{u}}^{n} \in \mathbb{R}^{N}$, where the coefficients on the interface $\underline{\mathsf{u}}_{\Gamma_{in}}^{n}$ satisfy the reduced Schur complement system 
\begin{equation}\label{eq:red_SC}
\underline{\Phi}_{n}^{T} \underline{A}_{SC} \underline{\Phi}_{n} \underline{\mathsf{u}}_{\Gamma_{in}}^{n} = \underline{\Phi}_{n}^{T} \underline{\mathsf{f}}_{SC} \qquad \text{of size } n \times n
\end{equation} 
and the DOFs of $\underline{\mathsf{u}}^{n}$ in the interior of $\Omega_{i}$, $i=1,2$ can be obtained in a standard manner via the definition of $\underline{A}_{SC}$.
The question of how to construct a reduced basis $\underline{\phi}_{1}, \hdots, \underline{\phi}_{n} \in \mathbb{R}^{N_{in}}$ which yields a rapidly convergent approximation and is even in some sense optimal will be addressed in the next section.

\section{Optimal port spaces}\label{sec:3}

Rather than assuming a priori knowledge about the shape of the global system associated with $\Omega_{gl}$, we wish to enable maximum flexibility in terms of system assembly on the user's side. In other words, we wish to supply the user with many components (or subdomains), each equipped with (local) reduced models, which the user can then use to build the desired system and thus implicitly define $\Omega_{gl}$. As a consequence, due to the a priori unknown geometry of $\Omega_{gl}$, we assume that the trace of the global solution $u_{gl}$ on $\Gamma_{out}$ is \emph{unknown} to us when constructing the reduced basis $\underline{\phi}_{1}, \hdots, \underline{\phi}_{n} \in \mathbb{R}^{N_{in}}$. We thus aim at approximating all local solutions of \eqref{eq:PDE} with \emph{arbitrary Dirichlet boundary conditions} on $\Gamma_{out}$. Before presenting the construction of the reduced basis in subsection \ref{subsec:opt_port_sp} we illustrate in a motivating example taken from \cite[Remark 3.3]{SmePat16} why we may hope to be able to find a low-dimensional port space that approximates the set of all local solutions well. 

\begin{remark}\label{rem:motivation}
We consider two components $\Omega_{i} \subset \mathbb{R}^{2}$, $i=1,2$ each of height $H$ in $x_2$ and length $L$ in $x_1$, such that $\Gamma_{out}$ is at $x_{1} = -L$ and $x_{1} = L$ and $\Gamma_{in}$ is at $x _{1}= 0$. We consider the Laplacian and impose homogeneous Neumann conditions on $x_2 = 0$ and $x_2 = H$ in both subdomains. Proceeding with separation of variables, we can infer that all local solutions of the PDE for this problem are of the form
\begin{equation}\label{beispiel}
u(x_{1},x_{2}) = a_{0} + b_{0}x_{1} + \sum_{n=1}^{\infty} \cos (n\pi\frac{x_{2}}{H}) \left[a_{n}\cosh(n\pi \frac{x_{1}}{H}) + b_{n}\sinh(n\pi \frac{x_{1}}{H})\right],
\end{equation}
where the coefficients $a_{n}, b_{n} \in \mathbb{R}$, $n=0,\hdots,\infty$ are determined by the Dirichlet data on $\Gamma_{out}$.
Thanks to the $\cosh$ function we can observe a very rapid and exponential decay of the local solutions \eqref{beispiel} in the interior of $\Omega$. Therefore, most of the local solutions \eqref{beispiel} have negligibly small values on $\Gamma_{in}$, which is why we expect a low-dimensional port space on $\Gamma_{in}$ to be able to provide a very good approximation of all local solutions \eqref{beispiel}. The construction procedure described below generalizes the separation of variables ansatz. 
\end{remark}

\subsection{Construction of optimal port spaces via a transfer operator}\label{subsec:opt_port_sp}

First, we address the case $f=0$; the general case will be dealt with at the end of this subsection. Motivated by the separation of variables procedure, and the fact that the global solution $u_{gl}$ on $\Omega_{gl}$ satisfies the PDE locally on $\Omega$, we consider the space of all local solutions of the PDE
\begin{equation}\label{eq:space of local solutions}
\mathcal{H}:= \{ w \, : \, \mathcal{A}w = 0 \enspace \text{in} \enspace \Omega, \enspace w = 0 \enspace \text{on} \enspace \Sigma_{D}\cap\partial\Omega\}.
\end{equation}
As in \cite{BabLip11,SmePat16,BuhSme18} we may then introduce a \emph{transfer operator} $\mathcal{T}: \mathcal{S} \rightarrow \mathcal{R}$ that takes arbitrary data $\zeta$ on $\Gamma_{out}$ as an input, solves the PDE $\mathcal{A}u = 0$ on $\Omega$ with that data $\zeta$ as Dirichlet boundary conditions on $\Gamma_{out}$, and finally restricts the local solution to $\Gamma_{in}$. Introducing the source and range spaces $\mathcal{S}:= \{ w|_{\Gamma_{out}} \, : \, w \in \mathcal{H}\}$ and $\mathcal{R}:= \{ (w - P_{\mathrm{ker}(\mathcal{A})}(w))|_{\Gamma_{in}} \, : \, w \in \mathcal{H}\}$ the transfer operator is thus defined as
\begin{equation}\label{eq:transfer_op_gamma}
\mathcal{T}(w|_{\Gamma_{out}})=\left(w-P_{\mathrm{ker}(\mathcal{A})}(w)\right)|_{\Gamma_{in}} \qquad \text{for} \enspace w \in \mathcal{H}.
\end{equation}
Here, $P_{\mathrm{ker}(\mathcal{A})}(w)$ denotes the orthogonal projection of $w$ on the kernel of the differential operator. Note that for instance for the Laplacian $\mathrm{ker}(\mathcal{A})$ equals the constant functions and in the case of linear elasticity $\mathrm{ker}(\mathcal{A})$ is the space of the rigid body motions. Following up Remark \ref{rem:motivation} note that the transfer operator allows us to assess how much of the data on $\Gamma_{out}$ reaches the inner interface $\Gamma_{in}$. It can then be shown that thanks to the Caccioppoli inequality $\mathcal{T}$ is compact and that certain eigenfunctions of $\mathcal{T}^{*}\mathcal{T}$ span the optimal port space, where $\mathcal{T}^{*}:\mathcal{R} \rightarrow \mathcal{S}$ denotes the adjoint operator (see \cite{BabLip11,SmePat16,Pinkus85} for details). Here, we use the concept of optimality in the sense of Kolmogorov \cite{Kol36}: A subspace $\mathcal{R}^{n} \subset \mathcal{R}$ of dimension at most $n$ for which holds
\begin{equation*}
d_{n}(\mathcal{T}(\mathcal{S});\mathcal{R}) = \sup_{\psi \in \mathcal{S}} \inf_{\zeta \in \mathcal{R}^{n}} \frac{\|\mathcal{T}\psi - \zeta\|_{\mathcal{R}}}{\|\psi\|_{\mathcal{S}}}
\end{equation*}
is called an optimal subspace for $d_{n}(\mathcal{\mathcal{T}}(\mathcal{S});\mathcal{R})$, where the Kolmogorov $n$-width $d_{n}(\mathcal{T}(\mathcal{S});\mathcal{R})$ is defined as 
\begin{equation*}
d_{n}(\mathcal{T}(\mathcal{S});\mathcal{R}) := \underset{\dim(\mathcal{R}^{n})=n}{\inf_{\mathcal{R}^{n}\subset \mathcal{R}}} \sup_{\psi \in \mathcal{S}} \inf_{\zeta \in \mathcal{R}^{n}} \frac{\|\mathcal{T}\psi - \zeta\|_{\mathcal{R}}}{\|\psi\|_{\mathcal{S}}}.
\end{equation*}
We summarize the findings about the optimal port spaces in the following proposition.

\begin{proposition}[Optimal port spaces \cite{SmePat16}]\label{theorem:eigenvalue problem}
The optimal port space is given by 
\begin{equation}\label{eq:optimal space}
\mathcal{R}^{n}:= \mathrm{span} \{\phi_{1}^{sp},...,\phi_{n}^{sp}\}, \qquad \text{where} \enspace \phi_{j}^{sp}=\mathcal{T}\varphi_{j},\quad j=1,...,n,
\end{equation}
and 
$\lambda_j$ are the largest $n$ eigenvalues and $\varphi_j$ the corresponding eigenfunctions
that satisfy the \emph{transfer eigenvalue problem}: Find $(\varphi_{j},\lambda_{j}) \in (\mathcal{S},\mathbb{R}^{+})$ such that
\begin{align}\label{eq:transfer eigenvalue problem}
(\,\mathcal{T}\varphi_{j}\,,\,\mathcal{T}w\,)_{\mathcal{R}} = \lambda_{j} (\,\varphi_{j}\,,\,w\,)_{\mathcal{S}} \quad \forall w \in \mathcal{S}.
\end{align}
Moreover, the following holds: 
\begin{equation}\label{eq:n-width_equals_eigenvalue}
d_{n}(\mathcal{T}(\mathcal{S});\mathcal{R}) = \sup_{\xi \in \mathcal{S}} \inf_{\zeta \in \mathcal{R}^{n}} \frac{\| \mathcal{T}\xi - \zeta\|_{\mathcal{R}}}{\|\xi\|_{\mathcal{S}}} = \sqrt{\lambda_{n+1}}.
\end{equation}
\end{proposition} 

\begin{remark}\label{separation of variables} 
We note that, as can be seen from \eqref{eq:transfer eigenvalue problem}, the optimal modes are those that maximize the energy on the inner interface $\Gamma_{in}$ relative to the energy they have on $\Gamma_{out}$. The optimal port space is thus spanned by the modes that relatively still contain the most information on $\Gamma_{in}$. For our motivating example discussed in Remark \ref{rem:motivation} we obtain $\mathcal{R}^{n}:= \mathrm{span}\{\cos(\pi\frac{x_{2}}{H}),\cos(2\pi\frac{x_{2}}{H}),\hdots,\cos(n\pi\frac{x_{2}}{H})\}$. Moreover, we can exploit the separation of variables solution to solve \eqref{eq:transfer eigenvalue problem} in closed form: $\lambda_{j} = \left(\cosh(L\sigma_{j-1})\right)^{-2}$, $j=1,2,3,...$, where the eigenproblem in $x_2$ in the separation of variables procedure yields separation constants $\sigma_{j}=(j\pi)/H$, $j=0,1,2,...$. This simple model problem also foreshadows the potentially very good performance of the associated optimal space \eqref{eq:optimal space} in light of \eqref{eq:n-width_equals_eigenvalue} and Proposition \ref{a priori bound}: we obtain exponential convergence.
\end{remark}

For $f\neq 0$ we solve the problem: Find $u^{f}$ such that 
\begin{equation*}
\mathcal{A}u^{f} = f \quad \text{in } \Omega \qquad \text{and} \qquad u = 0 \quad \text{on } \Gamma_{out}
\end{equation*}
and augment the space $\mathcal{R}^{n}$ with $u^{f}|_{\Gamma_{in}}$ to arrive at
\begin{equation}\label{eq:optimal space interface}
\mathcal{R}_{data,\mathrm{ker}}^{n}:=\mathrm{span} \{ \phi_{1}^{sp},...,\phi_{n}^{sp}, u^{f}|_{\Gamma_{in}},\eta_{1}|_{\Gamma_{in}},\hdots,\eta_{\dim(\mathrm{ker}(\mathcal{A}))}|_{\Gamma_{in}}\},
\end{equation}
where $\{ \eta_{1},\hdots,\eta_{\dim(\ker(\mathcal{A}))}\}$ denotes a basis for $\mathrm{ker}(\mathcal{A})$.

Using the optimal port space $\mathcal{R}_{data,\mathrm{ker}}^{n}$ within the static condensation procedure allows proving the following {\it a priori} error bound for the static condensation approximation. We note that Proposition \ref{a priori bound} gives a bound for the continuous analogon $u^{n}$ of $\mathsf{u}^{n}$, the latter being defined in subsection \ref{subsect:port red}. To simplify the notation we do not give a precise definition of $u^{n}$ and refer to that end to \cite{SmePat16}. Note however, that the convergence behavior of $\mathsf{u}^{n}$ towards $\mathsf{u}$ is very similar to the continuous setting, differing only due to the FE approximation. 

\begin{proposition}[{\it A priori} error bound \cite{SmePat16}]\label{a priori bound}
Let $u$ be the (exact) solution of \eqref{eq:PDE} and $u^{n}$ the continuous static condensation approximation employing the optimal port space $\mathcal{R}_{data,\mathrm{ker}}^{n}$. Moreover, denote with $\| \cdot \|_{\mathcal{E}}$ the norm induced by the bilinear form associated with the differential operator $\mathcal{A}$.
Then we have the following {\it a priori} error bound:
\begin{equation}\label{eq a priori bound}
\frac{\| u - u^{n} \|_{\mathcal{E}}}{\| u\|_{\mathcal{E}}} \leq C_{1}(\Omega) \,\sqrt{\lambda_{n+1}},
\end{equation}
where $\lambda_{n+1}$ is the $n+1$th eigenvalue of \eqref{eq:transfer eigenvalue problem} and $C_{1}(\Omega)$ is a constant which depends neither on $u$ nor on $u^{n}$.
\end{proposition}

\subsection{Approximation of the optimal spaces}

In this subsection we show how an approximation of the continuous optimal local spaces $\mathcal{R}_{data,\ker}^{n}$ can be computed with the FE method. First, in order to define a matrix form of the transfer operator we introduce DOF mappings $\underline{D}_{\Gamma_{out}\rightarrow \Omega} \in \mathbb{R}^{N \times N_{out}}$ and $\underline{D}_{\Omega\rightarrow \Gamma_{in}} \in \mathbb{R}^{N_{in}\times N}$ that map the DOFs on $\Gamma_{out}$ to the DOFs of $\Omega$ and the DOFs of $\Omega$ to the DOFs of $\Gamma_{in}$, respectively; $N_{out}$ denotes the number of DOFs on $\Gamma_{out}$. By denoting with $\underline{\zeta} \in \mathbb{R}^{N_{out}}$ the coefficients of a FE function $\zeta$ on $\Gamma_{out}$ and denoting by $\underline K_{\Omega}$ the matrix of the orthogonal projection $P_{ker(\mathcal{A}),\Omega}$ on $\ker(\mathcal{A})$ on $\Omega$ we obtain the following matrix representation $\underline{T} \in \mathbb{R}^{N_{in} \times N_{out}}$ of the transfer operator:
\begin{eqnarray}\label{eq:matrix_form_transfer_operator2}
\underline{T}\,\underline{\zeta} = 
\underline{D}_{\Omega\rightarrow \Gamma_{in}} \, 
\left(1 - \underline K_{\Omega} \right) \,
\underline{A}^{-1} \underline{D}_{\Gamma_{out}\rightarrow \Omega} \, \underline{\zeta}.
\end{eqnarray}
Finally, we denote by $\underline M_S$ the inner product matrix of the FE source space $S$ and by $\underline M_R$ the inner product matrix of the FE range space $R$. Possible inner products for $S$ and $R$ are the $L^{2}$-inner product and a lifting inner product. To obtain the latter we solve for instance for a function $\xi$ defined on $\Gamma_{in}$ the PDE on $\Omega_{1}$ and $\Omega_{2}$ numerically with Dirichlet data $\xi$ on $\Gamma_{in}$ and homogeneous Dirichlet boundary conditions on $\Gamma_{out}$ --- for further details we refer to \cite{SmePat16} and for the FE implementation to the Supplementary Materials of \cite{SmePat16}. The FE approximation of the transfer eigenvalue problem then reads as follows: Find the eigenvectors $\underline{\zeta}_{j} \in \mathbb{R}^{N_{S}}$ and the eigenvalues $\lambda_{j} \in \mathbb{R}^{+}_0$ such that
\begin{equation}\label{eq:matrix_version_transfer_eigenvalue_problem}
\underline{T}^{t}\underline M_{R}\underline{T} \,\underline{\zeta}_{j} = \lambda_{j}\, \underline{M}_{S} \,\underline{\zeta}_{j}.
\end{equation}
Note that in actual practice we would not assemble $\underline{T}$ but instead solve successively the linear system of equations
\begin{equation}\label{eq:local FE PDE}
\underline{A}\underline{u}_{i} = D_{\Gamma_{out}\rightarrow \Omega}\underline{e}_{i} \qquad \text{with the standard unit vectors } \underline{e}_{i}
\end{equation}
 and assemble $(\underline{T}^{T} \underline{M}_{R} \underline{T})_{i,j} = (\underline{D}_{\Omega\rightarrow \Gamma_{in}}\underline{u}_{j},\underline{D}_{\Omega\rightarrow \Gamma_{in}}\underline{u}_{i})_{R}.$ The coefficients of the FE approximation of the basis functions $\{\phi_{1}^{sp},...,\phi_{n}^{sp}\}$ of the optimal local approximation space $R^{n}:=\mathrm{span} \{\phi_{1}^{sp},...,\phi_{n}^{sp}\}$
are then given by
$
\underline{\phi}_{j}^{sp} = \underline{T}\,\underline{\zeta}_{j},$ $j =1,\hdots,n.
$
Adding the representation of the right-hand side and a basis of $\ker(\mathcal{A})$ yields the optimal space $R^{n}_{data,\mathrm{ker}}$.

\begin{remark}
We may also define the discrete transfer operator implicitly via \eqref{eq:matrix_form_transfer_operator2} and pass it together with its implicitly defined adjoint to a Lanczos method. This is in general much more favorable from a computational viewpoint compared to solving \eqref{eq:local FE PDE} $N_{out}$ times. However, it turns out that employing techniques from randomized linear algebra can be even more computationally beneficial than a Lanczos method as it requires only about $n$ local solutions of the PDE with random boundary conditions while yielding an approximation of the eigenvectors $\underline{\zeta}_{j}$ of \eqref{eq:matrix_version_transfer_eigenvalue_problem} at any required accuracy \cite{BuhSme18}. 
\end{remark}

\section{Extension to parameter-dependent problems and problems with geometric changes} \label{sec:4}
Many applications require a rapid simulation response for many different material parameters such as Young's modulus or a real-time simulation response for a different geometry such as a beam with a newly detected crack. Therefore, it is desirable to have a port-reduced static condensation procedure that is able to deal efficiently with parameter-dependent PDEs and geometric changes. Recall however that the optimal port space as presented in section \ref{sec:3} is based on the space of functions that solve the (now parametrized) PDE on a specific domain $\Omega$ and therefore also depends on the parameter and the geometry of $\Omega$. As constructing a new optimal port space ``from scratch'' for each new parameter value is in general not feasible, the goal of this section is to show how to construct a low-dimensional and quasi-optimal port space that is independent of the parameter and the geometry but yields an accurate approximation for the full parameter set and all geometries of interest. To that end, we present in subsection \ref{subsect:greedy} a spectral greedy approach which constructs a reduced basis to approximate the $n$ eigenspaces associated with the $n$ largest eigenvalues of the parameter and geometry dependent generalized (transfer) eigenvalue problem. Here, we slightly extend the spectral greedy algorithm introduced in \cite{SmePat16} to the case of varying geometries. At the beginning we state the parametrized PDE of interest in subsection \ref{subsect:parameter_problem_setting} and recall the port-reduced static condensation procedure for parameter-dependent PDEs in subsection \ref{subsect:sc_parameter}.

\subsection{Parametrized partial differential equations with geometric changes}\label{subsect:parameter_problem_setting}
We consider a setting where $\Omega(\Mu)$ accounts for different geometries such a beam, a beam with a crack or a beam with a hole. Note that in contrast to ``standard'' model order reduction appraoches we do accomodate here geometries that cannot be transformed in one another by a $C^{1}$-map. Moreover, we allow different discretizations in the interior of $\Omega_{i}(\Mu)$, $i=1,2$. However, we have to insist that the geometry of $\Gamma_{in}(\boldsymbol{\mu})$ is parameter-independent and that the meshes associated with all considered components coincide on $\Gamma_{in}(\boldsymbol{\mu})$; translation is of course possible resulting in the parameter dependency of $\Gamma_{in}(\Mu)$. 

In detail, we consider a discrete geometry parameter set $\mathcal{P}_{Geo}$ being the union of the considered geometries. Geometric changes via smooth maps can additionally be accounted for via the parameter-dependent operator $\mathcal{A}(\boldsymbol{\mu})$, where $\boldsymbol{\mu}$ belongs to the compact parameter set $\mathcal{P}_{PDE} \subset \mathbb{R}^{p}$. Again we assume that the port is not geometrically deformed. Then, we consider the following problem on $\Omega(\boldsymbol{\mu})$: For any $\boldsymbol{\mu} \in \mathcal{P}:=\mathcal{P}_{Geo}\times \mathcal{P}_{PDE}$ and given $f(\boldsymbol{\mu}) \in L^{2}(\Omega(\boldsymbol{\mu}))$ find $u(\boldsymbol{\mu}) $ such that
\begin{equation}\label{eq:PDE_para}
\mathcal{A}(\boldsymbol{\mu}) u(\boldsymbol{\mu})  = f(\boldsymbol{\mu})  \enspace \text{in } \Omega(\boldsymbol{\mu}), \quad \text{and} \quad u(\boldsymbol{\mu}) = u_{gl}(\boldsymbol{\mu}) \enspace \text{on } \Gamma_{out}(\boldsymbol{\mu}).
\end{equation}
Again, we introduce a conforming FE discretization to arrive at the linear system of equations $\underline{A}(\Mu) \underline{\mathsf{u}}(\Mu) = \underline{\mathsf{f}}(\Mu)$ of size $N(\Mu) \times N(\Mu)$ and FE approximation $\mathsf{u}(\Mu)$.

\subsection{Port reduced static condensation for parametrized equations}\label{subsect:sc_parameter}

We assume that we have given a \emph{parameter-independent} reduced port basis $\underline{\phi}_{1}, \hdots, \underline{\phi}_{m} \in \mathbb{R}^{N_{in}}$, $m \ll N_{in}$ that we store in the columns of the matrix $\underline{\Phi}_{m} \in \mathbb{R}^{N_{in} \times m}$. Proceeding as above we can then define a \emph{parameter-dependent port reduced static condensation approximation} \cite{EftPat13a} $\mathsf{u}^{m}(\Mu)$ with FE coefficients $\underline{\mathsf{u}}^{m}(\Mu) \in \mathbb{R}^{N}$, where the coefficients on the interface $\underline{\mathsf{u}}_{\Gamma_{in}}^{m}(\Mu)$ satisfy the parametrized and reduced Schur complement system 
\begin{equation}\label{eq:red_SC_para}
\underline{\Phi}_{m}^{T} \underline{A}_{SC}(\Mu) \underline{\Phi}_{m} \underline{\mathsf{u}}_{\Gamma_{in}}^{m}(\Mu) = \underline{\Phi}_{m}^{T} \underline{\mathsf{f}}_{SC}(\Mu) \qquad \text{of size } m \times m
\end{equation} 
and the Schur complement matrix $\underline{A}_{SC}(\Mu)\in \mathbb{R}^{N_{in} \times N_{in}}$ and the Schur complement right-hand $\underline{\mathsf{f}}_{SC}(\Mu)\in \mathbb{R}^{N_{in}}$ are defined as follows
\begin{align}
\nonumber \underline{A}_{SC}(\Mu) &= \underline{A}_{\Gamma_{in}}(\Mu) - \underline{A}_{\Gamma_{in},\Omega_1}^{T}(\Mu) \underline{A}_{\Omega_{1}}^{-1}(\Mu)  \underline{A}_{\Gamma_{in},\Omega_1}(\Mu) - \underline{A}_{\Gamma_{in},\Omega_2}^{T}(\Mu)\underline{A}_{\Omega_{2}}^{-1}(\Mu) \underline{A}_{\Gamma_{in},\Omega_2}(\Mu),\\[-1.5ex]
& \label{eq:SC_matrix_para}\\[-1.5ex]
\nonumber \underline{\mathsf{f}}_{SC}(\Mu) &= \underline{\mathsf{f}}_{\Gamma_{in}}(\Mu) - \underline{A}_{\Gamma_{in},\Omega_1}^{T}(\Mu) \underline{A}_{\Omega_{1}}^{-1}(\Mu)\underline{\mathsf{f}}_{\Omega_{1}}(\Mu) -  \underline{A}_{\Gamma_{in},\Omega_2}^{T}(\Mu)\underline{A}_{\Omega_{2}}^{-1}(\Mu)\underline{\mathsf{f}}_{\Omega_{2}}(\Mu).
\end{align}
We note that in order to facilitate a simulation response at low marginal cost one would in actual practice also use model order reduction techniques to approximate $\underline{A}_{\Omega_{i}}^{-1}(\Mu)  \underline{A}_{\Gamma_{in},\Omega_i}(\Mu)$, $i=1,2$. This is however not the topic of this paper and we refer for details to \cite{HuKnPa13}.

\subsection{Spectral greedy algorithm}\label{subsect:greedy}

The process defined in section \ref{sec:3} yields for every $\Mu \in \mathcal{P}$ the (optimal) port space $R^{n}_{data,\mathrm{ker}}(\Mu)$ for this specific parameter $\Mu \in \mathcal{P}=\mathcal{P}_{Geo}\times \mathcal{P}_{PDE}$. 
The spectral greedy algorithm as introduced in \cite{SmePat16} and which we extend here to the case of geometry changes constructs \emph{one} quasi-optimal parameter-independent port space $R^{m}$ which approximates those parameter-dependent spaces $R^{n}_{data,\mathrm{ker}}(\Mu)$ with a given accuracy on a finite dimensional training set $\Xi =\mathcal{P}_{Geo} \times \Xi_{PDE}$ with $\Xi_{PDE}  \subset \mathcal{P}_{PDE}$. In the spectral greedy algorithm we exploit the fact that, although the solutions on the component pair may vary significantly with the parameter $\Mu \in \mathcal{P}_{PDE}$ and the geometry, we expect that the port spaces $R^{n}_{data,\mathrm{ker}}(\Mu)$, and in particular the spectral modes that correspond to the largest eigenvalues, are much less affected by a variation in the parameter and changes in the geometry thanks to the expected very rapid decay of the higher eigenfunctions in the interior of $\Omega(\Mu)$. 

The spectral greedy as described in Algorithm \ref{algorithm spectral greedy} then proceeds as follows. After the initialization we compute for all $\Mu \in \Xi$ the parameter-dependent optimal port spaces $R^{n}_{data,\mathrm{ker}}(\Mu)$. Also in the parameter-dependent setting we can prove an {\it a priori} error bound \cite{SmePat16} for the error between $u(\Mu)$ and the continuous port-reduced static condensation approximation $u^{n}(\Mu)$ corresponding to the \emph{parameter-depedent} optimal port space $\mathcal{R}^{n}_{data,\ker}(\Mu)$:
\begin{equation}\label{eq:a priori para}
\frac{\| u(\Mu) - u^{n}(\Mu) \|_{\mathcal{E}(\Mu)} }{\|u(\Mu)\|_{\mathcal{E}(\Mu)}}\leq c_{1}(\Mu)c_{2}(\Mu) C_{1}(\Omega(\Mu),\Mu) \sqrt{\lambda_{n+1}(\Mu)}.
\end{equation}
Here, the norm $\| \cdot \|_{\mathcal{E}(\Mu)}$ is the norm induced by the parameter-dependent bilinear form associated with $\mathcal{A}(\Mu)$ and $c_{1}(\Mu)$ and $c_{2}(\Mu)$ are chosen such that we have $c_{1}(\Mu) \| \cdot \|_{\mathcal{E}(\bar{\Mu})} \leq \| \cdot \|_{\mathcal{E}(\Mu)} \leq c_{2}(\Mu) \| \cdot \|_{\mathcal{E}(\bar{\Mu})}$ for all $\Mu \in \mathcal{P}_{PDE}$ and a fixed reference parameter $\bar{\Mu} \in \mathcal{P}_{PDE}$.\footnote{We note that in order to prove \eqref{eq:a priori para} it is necessary to define the lifting inner product on the ports for one reference parameter $\bar{\Mu} \in \mathcal{P}_{Geo}$ and use the equivalence of the norm induced by the lifting inner product and the $H^{1/2}$-norm on the ports. Exploiting that the latter is the same for all considered geometries allows switching between the geometries in the proof.} To ensure that for every parameter $\Mu \in \Xi$ we include all necessary information that we need to obtain a good approximation for this specific parameter $\Mu$ we choose the dimension of $R^{n}_{data,\mathrm{ker}}(\Mu)$ for each $\Mu \in \Xi$ such that $c_{1}(\Mu)c_{2}(\Mu) C_{1}(\Omega(\Mu),\Mu) \sqrt{\lambda_{n+1}(\Mu)} \leq \frac{\varepsilon}{2}$ for a given tolerance $\varepsilon$. Although precise estimates for $C_{1}(\Omega(\Mu),\Mu)$ can be obtained, setting $C_{1}(\Omega(\Mu),\Mu)=1$ yields in general good results as another value would just result in rescaling $\varepsilon$. After having collected all vectors on $\Gamma_{in}$ that are essential to obtain a good approximation for all vectors $\underline{D}_{\Omega(\Mu)\rightarrow \Gamma_{in}}\underline{\mathsf{u}}(\Mu)$, $\Mu \in \Xi$, we must select a suitable basis from those vectors. This is realized in an iterative manner in Lines 5-14. 

\begin{algorithm}[t]
\caption{spectral greedy \cite{SmePat16}\label{algorithm spectral greedy}}
\SetKwData{Left}{left}\SetKwData{This}{this}\SetKwData{Up}{up}
\SetKwFunction{Union}{Union}\SetKwFunction{FindCompress}{FindCompress}
\SetKwInOut{Input}{input}\SetKwInOut{Output}{output}
\Input{train sample $\Xi \subset \mathcal{P}$, tolerance $\varepsilon$}
\Output{set of chosen parameters $P_{m}$, port space $R^m$}
\BlankLine
\textbf{Initialize} $P_{\dim(\mathrm{ker}(\mathcal{A}))}\leftarrow\emptyset, R^{\dim(\mathrm{ker}(\mathcal{A}))}\leftarrow\mathrm{span}\{\eta_{1}|_{\Gamma_{in}},\hdots,\eta_{\dim(\mathrm{ker}(\mathcal{A}))}|_{\Gamma_{in}}\}, m\leftarrow \dim(\mathrm{ker}(\mathcal{A}))$ \\
\ForEach{$\Mu \in \Xi$}{
Compute $R^{n}_{data,\mathrm{ker}}(\Mu)$ such that $c_{1}(\Mu)c_{2}(\Mu) C_{1}(\Omega(\Mu),\Mu)\sqrt{\lambda_{n+1}(\Mu)}\leq \frac{\varepsilon}{2}$.\label{alg compute Lambda}
}
\BlankLine
\While{true}{
\If{$\max_{\Mu \in \Xi} E(S(R^{n}_{data,\mathrm{ker}}(\Mu)),R^m) \leq \varepsilon/(\varepsilon + 2C_{2}(\Omega(\Mu),\Mu)c_{1}(\Mu)c_{2}(\Mu))$ \label{alg thres}}{
\textbf{return}\\
}
$\Mu^{*}\leftarrow \arg \max_{\Mu \in \Xi} E(S(R^{n}_{data,\mathrm{ker}}(\Mu)),R^m)$\label{alg mu max}\\
$P_{m+1} \leftarrow P_{m} \cup \Mu^{*}$\\
$\kappa \leftarrow \arg \sup_{\rho \in S(R^{n}_{data,\mathrm{ker}}(\Mu^{*}))} \inf_{\zeta \in R^m} \|\rho - \zeta\|_{R}$\label{kappa}\\
$R^{m+1} \leftarrow R^m + \mathrm{span}\{\kappa\}$ \label{alg extend}\\
$ m \leftarrow m+1$}
\end{algorithm}

In each iteration we first identify in Line \ref{alg mu max} the port space $R^{n}_{data,\mathrm{ker}}(\Mu^{*})$ that maximizes the deviation 
\begin{equation*}
E(S(R^{n}_{data,\mathrm{ker}}(\Mu)),R^m):= \sup_{\xi \in S(R^{n}_{data,\mathrm{ker}}(\Mu))} \inf_{\zeta \in R^{m}} \|\xi - \zeta\|_{R}, \quad \Mu \in \Xi,
\end{equation*}
where possible choices of $S(R^{n}_{data,\mathrm{ker}}(\Mu))\subset R^{n}_{data,\mathrm{ker}}(\Mu)$ will be discussed below. Subsequently, we determine in Line \ref{kappa} the function $\kappa \in S(R^{n}_{data,\mathrm{ker}}(\Mu^{*}))$ that is worst approximated by the space $R^{m}$ and enhance $R^{m}$ with the span of $\kappa$. The spectral greedy algorithm terminates if for all $\Mu \in \Xi$ we have 
\begin{equation}\label{eq:stop_greedy}
\max_{\Mu \in \Xi} E(S(R^{n}_{data,\mathrm{ker}}(\Mu)),R^{m}) \leq \varepsilon/(\varepsilon + 2C_{2}(\Omega(\Mu),\Mu)c_{1}(\Mu)c_{2}(\Mu))
\end{equation}
for a constant $C_{2}(\Omega(\Mu),\Mu)$, which can in general be chosen equal to one. A slight modification of the stopping criterion \eqref{eq:stop_greedy} and a different scaling of $\varepsilon$ in the threshold for the a priori error bound in Line \ref{alg compute Lambda} allows to prove that after termination of the spectral greedy we have \cite{SmePat16}
\begin{equation}
\|u(\Mu) - u^{m}(\Mu)\|_{\mathcal{E}(\Mu)}/\|u(\Mu) \|_{\mathcal{E}(\Mu)} \leq \varepsilon,
\end{equation}
where $u^{m}(\Mu)$ is the continuous port-reduced static condensation approximation corresponding to $\mathcal{R}^{m}$; $\mathcal{R}^{m}$ being the continuous outcome of the spectral greedy. 

\paragraph*{Choice of the subset $S(R^{n}_{data,\mathrm{ker}}(\Mu))$}

First, we emphasize that in contrast to the standard greedy as introduced in \cite{VePrRoPa03} we have an ordering of the basis functions in $R^{n}_{data,\mathrm{ker}}(\Mu)$ in terms of their approximation properties thanks to the transfer eigenvalue problem. To obtain a parameter-independent port space that yields a (very) good static condensation approximation already for moderate $m$ it is therefore desirable that the spectral greedy algorithm selects the more important basis functions sooner rather than later during the \texttt{while}-loop. The sorting of the basis functions in terms of their approximation properties is implicitly saved in their norms as $\|\phi_{j}(\Mu)\|_{R}^{2} = \lambda_{j}(\Mu)
$, $j=1,\hdots,n$ where $\phi_{j}(\Mu)$ denotes the spectral basis of $R^{n}_{data,\mathrm{ker}}(\Mu)$. As suggested in \cite{SmePat16} we thus propose to consider 
\begin{equation}\label{choice of subset in greedy}
S(R^{n}_{data,\mathrm{ker}}(\Mu)):= \{ \zeta(\Mu) \in R^{n}_{data,\mathrm{ker}}(\Mu)\,:\,  \sum_{i=1}^{\dim(R^{n}_{data,\mathrm{ker}}(\Mu))} (\underline{\zeta}(\Mu)_{i})^{2} \leq 1\} 
\end{equation}
with $\zeta(\Mu) = \sum_{i=1}^{\dim(R^{n}_{data,\mathrm{ker}}(\Mu))} \underline{\zeta}(\Mu)_{i} \phi_{i}(\Mu)$. The deviation $E(S(R^{n}_{data,\mathrm{ker}}(\Mu)),R^m)$ can then be computed by solving the eigenvalue problem: Find $(\underline{\psi}_{j}(\Mu),\sigma_{j}(\Mu)) \in (\mathbb{R}^{\dim(R^{n}_{data,\mathrm{ker}}(\Mu))},\mathbb{R}^{+})$ such that
\begin{align}\label{eigenvalue problem matrix}
\underline{Z}(\Mu)\underline{\psi}_{j}(\Mu) &= \sigma_{j}(\Mu) \underline{\psi}_{j}(\Mu),\\
\text{where} \enspace \underline{Z}_{i,l}(\Mu)&:= (\phi_{l}(\Mu) - \sum_{k=1}^{m}(\phi_{l}(\Mu),\phi_{k})_{R}\phi_{k},\phi_{i}(\Mu) - \sum_{k=1}^{m}(\phi_{i}(\Mu),\phi_{k})_{R}\phi_{k})_{R},
\end{align}
where $\phi_{k}$ denotes the basis of $R^{m}$ and the underscore denotes the coefficients of a vector in $R^{n}_{data,\mathrm{ker}}(\Mu)$ expressed in the spectral basis $\phi_{l}(\Mu)$. We thus obtain 
$
E(S(R^{n}_{data,\mathrm{ker}}(\Mu)),R^m) = \sqrt{\sigma_{1}(\Mu)}
$, for all $\Mu \in \Xi$,
and $\kappa = \psi_{1}(\Mu^{*})$ at each iteration. To further motivate this choice of $S(R^{n}_{data,\mathrm{ker}}(\Mu))$ let us assume that all spectral modes in $R^{n}_{data,\mathrm{ker}}(\Mu)$ are orthogonal to the space $R^m$ for all $\Mu \in \Xi$, which is the case for instance for $m=\dim(R^{n}_{data,\mathrm{ker}}(\Mu))$ but also often for higher $m$. In this case the matrices $\underline{Z}(\Mu)$ reduce to diagonal matrices with diagonal entries $\underline{Z}_{i,i}(\Mu)=\|\phi_{i}(\Mu)\|_{R}^{2}$, $i=1,...,\dim(R^{n}_{data,\mathrm{ker}}(\Mu))$, $\Mu \in \Xi$. A spectral greedy based on $E(S(R^{n}_{data,\mathrm{ker}}(\Mu)),R^m)$ would therefore select the parameter $\Mu^{*}$ such that the associated function $\psi_{1}(\Mu^{*})$ has maximal energy with respect to the $(\cdot , \cdot )_{R}$-inner product. Note that this is consistent with our aim to include the weighting induced by the transfer eigenvalue problem into the basis selection process by the spectral greedy. 

\begin{remark}
Note that were we to consider the norm $\| \cdot \|_{R}$ in \eqref{choice of subset in greedy} the sorting of the spectral basis $\phi_{i}(\Mu)$ of $R^{n}_{data,\mathrm{ker}}(\Mu)$ in terms of approximation properties is neglected in the \texttt{while} loop of Algorithm \ref{algorithm spectral greedy}; for further explanations see \cite{SmePat16}. As a consequence it may and often would happen in actual practice, also due to numerical inaccuracies, that a spectral greedy algorithm based on the $\| \cdot \|_{R}$ norm in \eqref{choice of subset in greedy} selects first functions that have been marked by the transfer eigenvalue problem as \emph{less} important. Therefore, we would observe an approximation behaviour of the static condensation approximation based on the so constructed port space that is not satisfactory for moderate $m$. 
\end{remark}

\section{Approximating the whole system associated with $\Omega_{gl}$}\label{sec:5}
To allow a maximal topological flexibility during assembly of the system associated with $\Omega_{gl}$, we assume that we neither know the size, the composition, nor the shape of the system when generating the reduced model. Therefore, we perform the spectral greedy algorithm for all interfaces that may appear in the large structure on the component pairs that share the interface. Multiplying the left-hand side of the inequality in Line \ref{alg compute Lambda} in Alg.~\ref{algorithm spectral greedy} and 
$2C_{2}(\Omega(\Mu),\Mu)c_{1}(\Mu)c_{2}(\Mu)$ in Line \ref{alg thres} in Alg.~\ref{algorithm spectral greedy} by an estimate for the number of times we expect the interface to appear in the large system ensures that the relative approximation error on the whole domain $\Omega_{gl}$ associated with the system will lie below $\varepsilon$ (see \cite{SmePat16} for the proof).\footnote{As indicated above it is necessary to slightly modify the spectral greedy algorithm to prove convergence.} We note that in actual practice numerical experiments show a very weak scaling in the number of ports such that the scaling might not be necessary \cite{SmePat16}.

\section{{\it A posteriori} error estimation}\label{sec:6}

In order to assess after the detection of a new defect in the assembled system whether the quality of the reduced port space is still sufficient we wish to have an {\it a posteriori} error estimator for the error between the port reduced solution $\mathsf{u}^{m}(\Mu)$ and the FE solution $\mathsf{u}(\Mu)$ available. To that end, we employ the error estimator derived in \cite{Sme14}. We exploit that the FE solution satisfies a \emph{weak flux continuity} at the interface $\Gamma_{in}$ 
\begin{equation}\label{eq:weak flux cont FE}
\underline{\mathsf{f}}_{SC}(\Mu) - \underline{A}_{SC}(\Mu) \underline{\mathsf{u}}_{\Gamma_{in}}(\Mu)=0.
\end{equation}
Regarding the term ``weak flux continuity'' we recall first that \eqref{eq:weak flux cont FE} is the discrete version of a Steklov-Poincar\'{e} interface equation. The latter is the weak counterpart of the Neumann condition $\frac{\partial u|_{\Omega_{1}}}{\partial n}=\frac{\partial u|_{\Omega_{2}}}{\partial n}$ on $\Gamma_{in}$ for the outer normal $n$, requiring continuity of the flux across the interface. For further details we refer to \cite{QuaVal05}.

Also the reduced solution $\underline{\mathsf{u}}^{m}(\Mu)$ satisfies a weak flux continuity with respect to the reduced test space:
 \begin{equation*}
\Phi_{m}^{T}\underline{\mathsf{f}}_{SC}(\Mu) - \Phi_{m}^{T}\underline{A}_{SC}(\Mu) \Phi_{m}\underline{\mathsf{u}}^{m}_{\Gamma_{in}}(\Mu)=0.
\end{equation*}
However, the reduced solution $\underline{\mathsf{u}}^{m}(\Mu)$ does \emph{not} satisfy a weak flux continuity with respect to the full test space:
 \begin{equation}\label{eq:no weak flux}
\underline{\Phi}_{N_{in}}^{T}\underline{\mathsf{f}}_{SC}(\Mu) - \underline{\Phi}_{N_{in}}^{T}\underline{A}_{SC}(\Mu) \Phi_{m}\underline{\mathsf{u}}^{m}_{\Gamma_{in}}(\Mu)\neq 0.
\end{equation}
Here, the first $m$ columns of $\underline{\Phi}_{N_{in}} \in \mathbb{R}^{N_{in} \times N_{in}}$ contain the basis $\underline{\phi}_{1},\hdots,\underline{\phi_{m}}$ generated by the spectral greedy and the remainder spans the orthogonal complement of $R^{m}$. Note that the left-hand side in \eqref{eq:no weak flux} can also be interpreted as a \emph{residual} on $\Gamma_{in}$. We use the violation of the weak flux continuity in \eqref{eq:no weak flux} to assess how much the reduced solution differs from the FE solution at the interface $\Gamma_{in}$. To utilize this information for {\it a posteriori} error estimation in \cite{Sme14} the concept of conservative fluxes defined according to Hughes et al. \cite{HuEnMaLa00} is adapted to the setting of port reduction. In a slight generalization of \cite{Sme14} we define the \emph{jump of the conservative flux} $\underline{\zeta}^{m}(\mu)$ as the solution of
\begin{equation}\label{eq:jump_cons_flux}
\underline{\Phi}_{N_{in}}^{T} \underline{M}_{R} \underline{\Phi}_{N_{in}} \underline{\zeta}^{m}(\mu) = \underline{\Phi}_{N_{in}}^{T}\underline{\mathsf{f}}_{SC}(\Mu) - \underline{\Phi}_{N_{in}}^{T}\underline{A}_{SC}(\Mu) \Phi_{m}\underline{\mathsf{u}}^{m}_{\Gamma_{in}}(\Mu).
\end{equation}
If $\phi_{1},\hdots,\phi_{m}$ are orthonormal with respect to the inner product in $R$ the linear system of equations \eqref{eq:jump_cons_flux} simplifies to 
\begin{equation}\label{eq:jump_cons_flux_simple}
 \underline{\zeta}^{m}(\mu) = \underline{\Phi}_{N_{in}}^{T}\underline{\mathsf{f}}_{SC}(\Mu) - \underline{\Phi}_{N_{in}}^{T}\underline{A}_{SC}(\Mu) \Phi_{m}\underline{\mathsf{u}}^{m}_{\Gamma_{in}}(\Mu).
\end{equation}
The computation of the jump of the conservative flux thus reduces to assembling the residual. Therefore, the computational costs scale linearly in $(N_{in} - m)$ and $m$.
\begin{proposition}[A posteriori error estimator for port reduction \cite{Sme14}]
Equip $R$ with the $L^{2}$-norm and define
\begin{equation}\label{def delta}
\Delta^{m}(\mu) := \frac{(\max_{i=1,2}c_{t^{*},i})\sqrt{1 +c_{p}^{2}}}{\alpha_{app}(\mu)}\|\zeta^{m} \|_{L^{2}(\Gamma_{in})}, 
\end{equation}
where $c_{t^{*},i}$ is the discrete trace constant in $\|v\|_{L^{2}(\Gamma_{in})} \leq c_{t^{*},i} \|v \|_{H^{1}(\Omega_{i})}$, $c_{p}$ is the constant in the Poincar\'{e}-Friedrichs-inequality, and $\alpha_{app}(\mu)$ an approximation of the FE coercivity constant $\alpha_{h}(\mu)$ of the bilinear form associated with $\mathcal{A}(\Mu)$. If $\alpha_{app}(\mu) \leq \alpha_{h}(\mu)$, there holds
\begin{equation}\label{delta}
\frac{1}{\gamma_{h}(\mu)c_{h}h^{-1/2}c_{a}}\enspace \Delta^{m}(\mu)\leq \| \nabla \mathsf{u}(\Mu) - \nabla \mathsf{u}(\Mu) \|_{L^{2}(\Omega_{\Mu})} \leq \Delta^{m}(\mu),
\end{equation}
where $c_{a}$ is the continuity constant of the discrete extension operator, $c_{h}$ is the constant in the inverse inequality $\|v\|_{H^{1/2}(\Gamma_{in})} \leq c_{h}h^{-1/2}\|v\|_{L^{2}(\Gamma_{in})}$, and $\gamma_{h}(\mu)$ the FE continuity constant of the bilinear form associated with $\mathcal{A}(\Mu)$.
\end{proposition}
\begin{remark}[Error estimation on $\Omega_{gl}$]
Let us assume that the system associated with $\Omega_{gl}$ has $P^{\Gamma}$ ports, which are denoted by $\Gamma_{p}$, $p=1,\hdots,P^{\Gamma}$. Moreover, denote by $\zeta_p^{m}$ the jump of the conservative flux at port $\Gamma_{p}$.
Then we can define an error estimator on $\Omega_{gl}$ as follows \cite{Sme14}:
\begin{equation}\label{eq:delta_2}
\Delta^{m}(\Mu) := \frac{c_{t^{*}} \sqrt{1 +c_{p}^{2}}}{\alpha_{app}(\mu)}\left( \sum_{p=1}^{P^{\Gamma}}\|\zeta_p^{m} \|_{L^{2}(\Gamma_{p})}^{2}\right)^{1/2}.
\end{equation}
Here, $c_{t^{*}}$ denotes the maximum over the discrete trace constants in all components, where we estimate the $L^{2}$-norm on all ports of that compontent against the $H^{1}$-norm on that component. $c_{p}$ is the constant in the Poincar\'{e}-Friedrichs-inequality with respect to $\Omega_{gl}$. Again, one can show that the effectivity of the error estimator \eqref{eq:delta_2} is bounded \cite{Sme14}. Moreover, we have that the effectivity of all local error indicators defined as in \eqref{def delta} is bounded. Those local error indicators associated with one port in the system can thus be used within an adaptive scheme to decide where to enrich the port space first.
\end{remark}

We note that due to coercivity constant and the constant $c_{p}$ the effectivities of $\Delta^{m}(\Mu)$ are in general rather high. However, in \cite{PatSme15} an error estimator is presented, which is solely based on local constants and in consequence provides a very sharp bound for the error. We finally note that the {\it a posteriori} error estimator introduced in \cite{Sme14} also assess the error due to an RB approximation of $\underline{A}_{\Omega_{i}}^{-1}(\Mu)  \underline{A}_{\Gamma_{in},\Omega_i}(\Mu)$, $i=1,2$ in \eqref{eq:red_SC_para}.

\section{Numerical Experiments}\label{sec:7}

In this section we investigate the performance of the optimal port space $R^{n}_{data,\mathrm{ker}}(\Mu)$ for changing geometries as occurring in structural health monitoring. We demonstrate in subsection \ref{subsect:num reuse} that we can use a port space generated on a component pair of two un-defecive (I-)beams also for a component pair with a defect such as a crack, obtaining a relative approximation error of less than $10^{-3}$. Subsequently, we investigate the performance of the spectral greedy algorithm for geometry changes in subsection \ref{subsect:num_greedy}. We begin in subsection \ref{subsect:benchmark} with the description of our benchmark problem: isotropic, homogeneous linear elasticity.

For the implementation we used the finite element library \texttt{libMesh} \cite{Kiretal06} including \texttt{rbOOmit} \cite{KnePet11}. The eigenvalue problems in the transfer eigenvalue problem and the computation of the deviation have been computed with the \texttt{Eigen} library \cite{Eigen}.

\subsection{Benchmark problem: Isotropic, homogeneous linear elasticity}\label{subsect:benchmark}

We assume that $\Omega(\Mu) \in \mathbb{R}^{d}$, $d=2,3$, $\bar{\Omega}(\Mu)=\bar{\Omega}_{1} \cup \bar{\Omega}_{2}(\Mu)$ is filled with an isotropic, homogeneous material and consider defects in the sense that $\Omega_{2}(\Mu)$ may have say a hole or a crack with a boundary $\Gamma_{defect}(\Mu)\subset\Gamma_{N}(\Mu)$. We consider the following linear elastic boundary value problem: Find the displacement vector $u(\Mu)$ and the Cauchy stress tensor $\sigma(u(\Mu))$ such that
\begin{align}
\nonumber - \nabla\cdot \sigma(u(\Mu)) &= 0 \quad &\text{in} \enspace \Omega(\Mu),\\
\label{strong form}\sigma(u(\Mu)) \cdot n(\Mu) &= 0  \quad &\text{on} \enspace \Gamma_{N}(\Mu),\\
\nonumber u(\Mu) &= g(\Mu) \quad &\text{on} \enspace \Gamma_{D},
\end{align}
where $g$ is a given Dirichlet boundary condition on the displacement. 

Thanks to Hooke's law we can express for a linear elastic material the Cauchy stress tensor as $\sigma(u(\Mu)) = C : \varepsilon(u(\Mu))$, where $C$ is the stiffness tensor, $\varepsilon(u(\Mu)) = 0.5 (\nabla u(\Mu) + (\nabla u(\Mu))^{T})$ is the infinitesimal strain tensor, and the colon operator $:$ is defined as $C : \varepsilon(u(\Mu)) = \sum_{i,j=1}^{2}C_{ij}\varepsilon_{ij}(u(\Mu))$. We assume in two spatial dimensions, i.e. for $d=2$, that the considered isotropic, homogeneous material is under plane stress. Therefore, the stiffness tensor can be written as 
\begin{align*}
C_{ijkl} = 
\begin{cases}
\frac{\nu}{(1 - \nu)^{2}}\delta_{ij}\delta_{kl} + \frac{1}{2(\nu + 1)}(\delta_{ik}\delta_{jl} + \delta_{il}\delta_{jk}), \enspace 1 \leq i,j,k,l \leq 2, \enspace &\text{if } d = 2,\\
\frac{\nu}{(1+\nu)(1 -2\nu)}\delta_{ij}\delta_{kl} + \frac{1}{2(1 + \nu)}(\delta_{ik}\delta_{jl} + \delta_{il}\delta_{jk}), \enspace 1 \leq i,j,k,l \leq 3, \enspace &\text{if } d = 3,
\end{cases}
\end{align*}
where $\delta_{ij}$ denotes the Kronecker delta and we choose Poisson's ratio $\nu =0.3$. We only consider parameters due to geometry changes such as a replacement of a beam with a cracked beam and no material parameters; therefore we have $\mathcal{P}_{PDE}=\emptyset$. As indicated in subsection \ref{subsect:parameter_problem_setting} we discretize the weak form of \eqref{strong form} by a conforming FE discretization. 

The kernel of $\mathcal{A}$ for the present example equals the three-dimensional space of rigid body motions for $d=2$ and the six-dimensional space of rigid body motions for $d=3$. To construct a port space $R^{n}_{data,\mathrm{ker}}(\Mu)$ on $\Gamma_{in}(\Mu)$ for each parameter we follow the procedure described in section \ref{sec:3}, where we use a lifting inner product (for further details on the latter see \cite{SmePat16}). As we do not consider a load here, we obtain $\dim(R^{n}_{data,\mathrm{ker}})=n+3$ for $d=2$ and $\dim(R^{n}_{data,\mathrm{ker}})=n+6$ for $d=3$. In order to construct one joint port space on $\Gamma_{in}(\Mu)$ we use the spectral greedy algorithm described in subsection \ref{subsect:greedy} using the $L^{2}$-inner product on $\Gamma_{in}$ and $C_{1}(\Omega(\Mu),\Mu)=C_{2}(\Omega(\Mu),\Mu)=1$; note that we have $c_{1}(\Mu)=c_{2}(\Mu)$ thanks to $\mathcal{P}_{PDE}=\emptyset$.

\subsection{Reusing the port space for a component with different geometry}\label{subsect:num reuse}
\begin{figure}[t]
\begin{center}
\includegraphics[scale = 0.33]{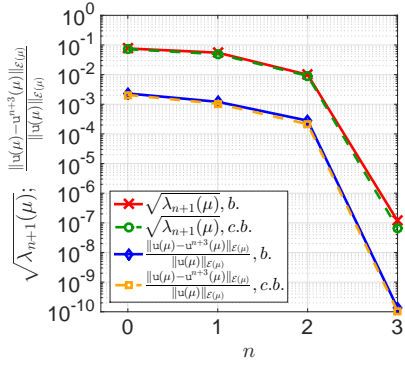}
\caption{Eigenvalues $\sqrt{\lambda_{n+1}(\Mu)}$ for the beam (b.) and the cracked beam (c. b.) and the average relative error $\| \mathsf{u}(\Mu) - \mathsf{u}^{n+3}(\Mu)\|_{\mathcal{E}(\Mu)}/\|\mathsf{u}(\Mu)\|_{\mathcal{E}(\Mu)}$ if the respective spectral modes are employed in the static condensation approximation.\label{fig eigenvalues}}
\end{center}
\end{figure}

To provide a simulation response at low marginal cost it would be desirable if we could reuse the port space generated for certain geometries also for other geometries. Therefore, we investigate here the effect on the relative approximation error if we construct a port space on a component pair of two un-defective beams and use that port space for the approximation on a component pair consisting of one un-defective beam on $\Omega_{1}$ and various defective beams on $\Omega_{2}(\Mu)$. 

\begin{figure}[t]
\begin{minipage}{0.5\textwidth}
\begin{center}
\subfloat[{\footnotesize cracked beam}]{\includegraphics[scale = 0.20]{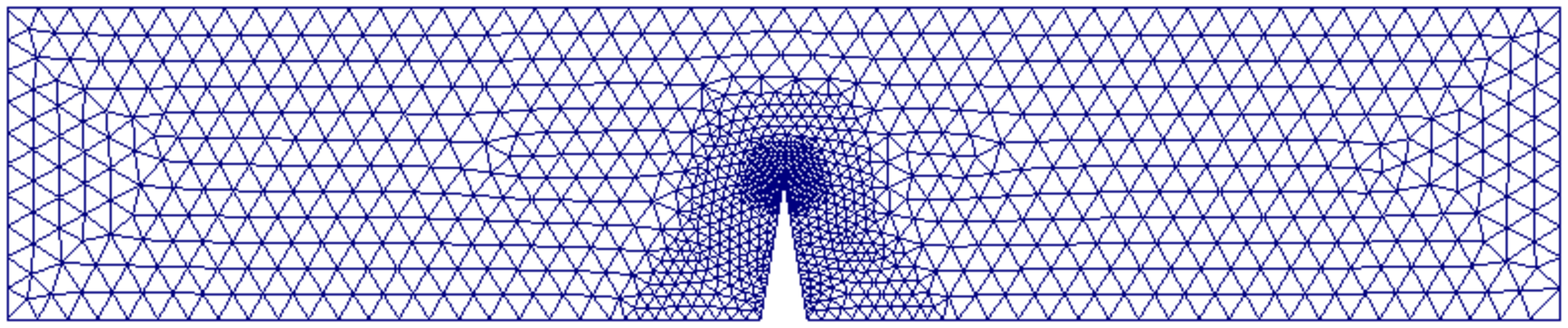} \label{fig_no_cut_out}}\quad
\subfloat[{\footnotesize beam with shifted crack}]{\includegraphics[scale = 0.2]{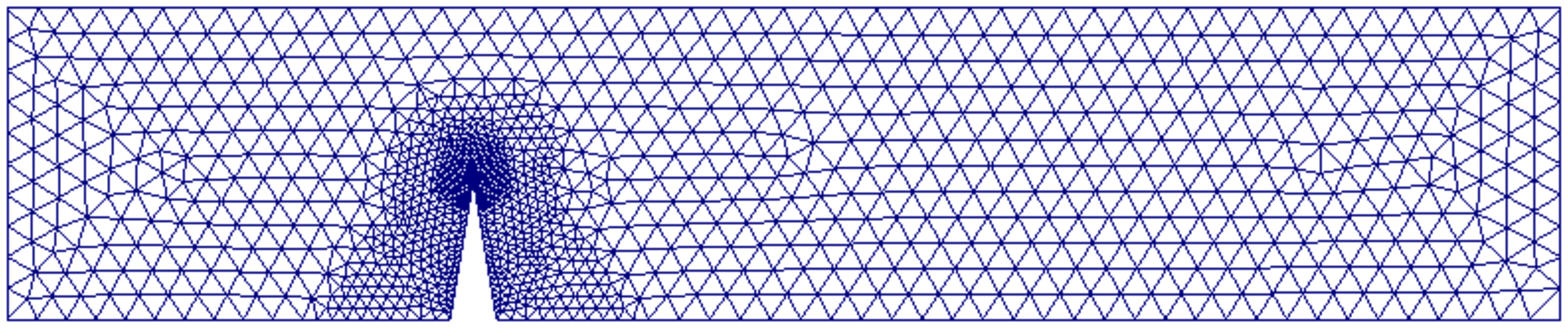} \label{fig_hex_mesh}}\quad
\subfloat[{\footnotesize beam with hole}]{\includegraphics[scale = 0.20]{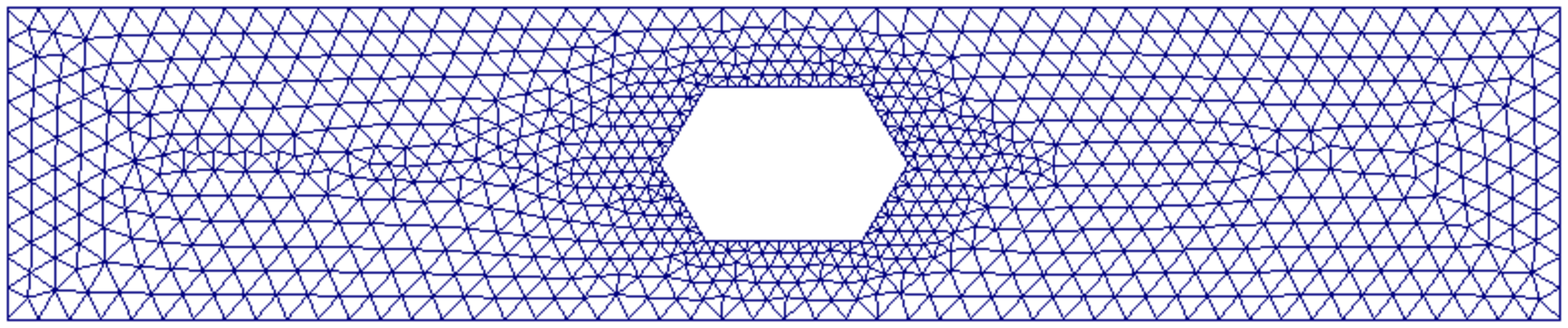} \label{fig_slit_mesh}} \quad
\caption{{\footnotesize Different component meshes}\label{fig_meshes}}
\end{center}
\end{minipage}
\quad
\begin{minipage}{0.45\textwidth}
\begin{center}
\includegraphics[scale = 0.33]{rel_error.eps}
\caption{{\footnotesize Average relative error $\| \mathsf{u}(\Mu) - \mathsf{u}^{n+3}(\Mu)\|_{\mathcal{E}(\Mu)}/\|\mathsf{u}(\Mu)\|_{\mathcal{E}(\Mu)}$ for various defective components.}\label{fig rel error}}
\end{center}
\end{minipage}
\end{figure}

To this end, we first connect components associated with the subdomains  $\Omega_{1}= (-7.5,-2.5)\times (-0.5,0.5)$ and $\Omega_{2}= (-2.5,2.5)\times (-0.5,0.5)$, i.e. two un-defective beams, and construct the associated port space. In the online stage we then prescribe random Dirichlet boundary conditions drawn uniformly from the interval $[-5,5]$ on the non-shared ports and verify that the average relative error $\| \mathsf{u}(\Mu) - \mathsf{u}^{n+3}(\Mu) \|_{\mathcal{E}(\Mu)}/\| \mathsf{u}(\Mu) \|_{\mathcal{E}(\Mu)}$ over 20 realizations exhibits nearly the same convergence behavior as $\sqrt{\lambda_{n+1}(\Mu)}$ (see Fig.~\ref{fig eigenvalues}). Subsequently, we replace the component associated with $\Omega_{2}$ by defective components: a cracked beam, a beam where the crack is shifted towards the shared port, and a beam with a hole; the corresponding component meshes are depicted in Fig.~\ref{fig_meshes}. Again we prescribe random Dirichlet boundary conditions on the non-shared ports and analyze the behavior of the average relative error for an increasing number of spectral modes that have been constructed by connecting two un-defective components. As anticipated the convergence behavior of the static condensation approximation for the defective components is (much) worse as that of the un-defective component (see Fig.~\ref{fig rel error}). Analyzing the convergence behavior of the static condensation approximation for the cracked beam using a port space that has been constructed by connecting a beam with a cracked beam (see Fig.~\ref{fig eigenvalues}) demonstrates that this worse convergence behavior is solely due to the fact that for the results in Fig.~\ref{fig rel error} we have employed the port space for the un-defective components. However, we emphasize that already for six spectral modes (including the three rigid body modes) we obtain for the defective components a relative error of about $10^{-5}$ (see Fig.~\ref{fig rel error}), which is very satisfactory. Moreover, we observe that the error increases only slightly when we shift the crack towards the shared port. 

Therefore, we conclude that in two space dimensions reusing the port space of the un-defective component yields a sufficiently accurate static condensation approximation. However, it should be noted that the port space contains only six port modes and is therefore rather small.  

\begin{figure}[t]
\begin{minipage}{0.6\textwidth}
\subfloat[{\footnotesize I-Beam}]{\includegraphics[scale = 0.350]{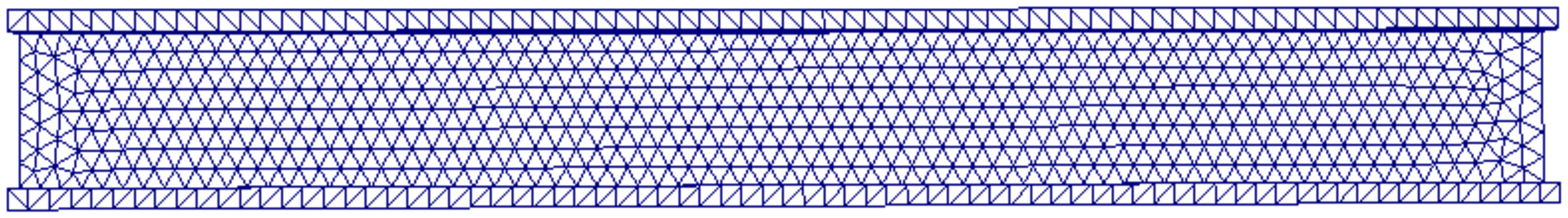} \label{fig_ibeam}}\\
\subfloat[{\footnotesize cracked I-Beam}]{\includegraphics[scale = 0.35]{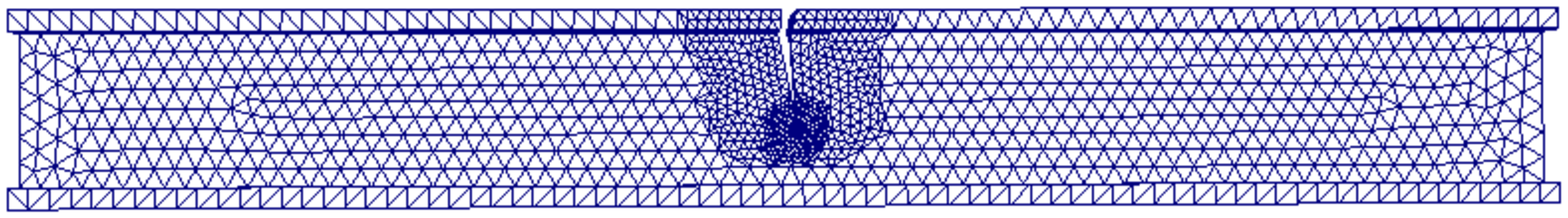} \label{fig_cracked_ibeam}}
\end{minipage}
\qquad \qquad\qquad
\begin{minipage}{0.3\textwidth}
\subfloat[{\footnotesize mesh on $\Gamma_{in}$}]{\includegraphics[scale = 0.20]{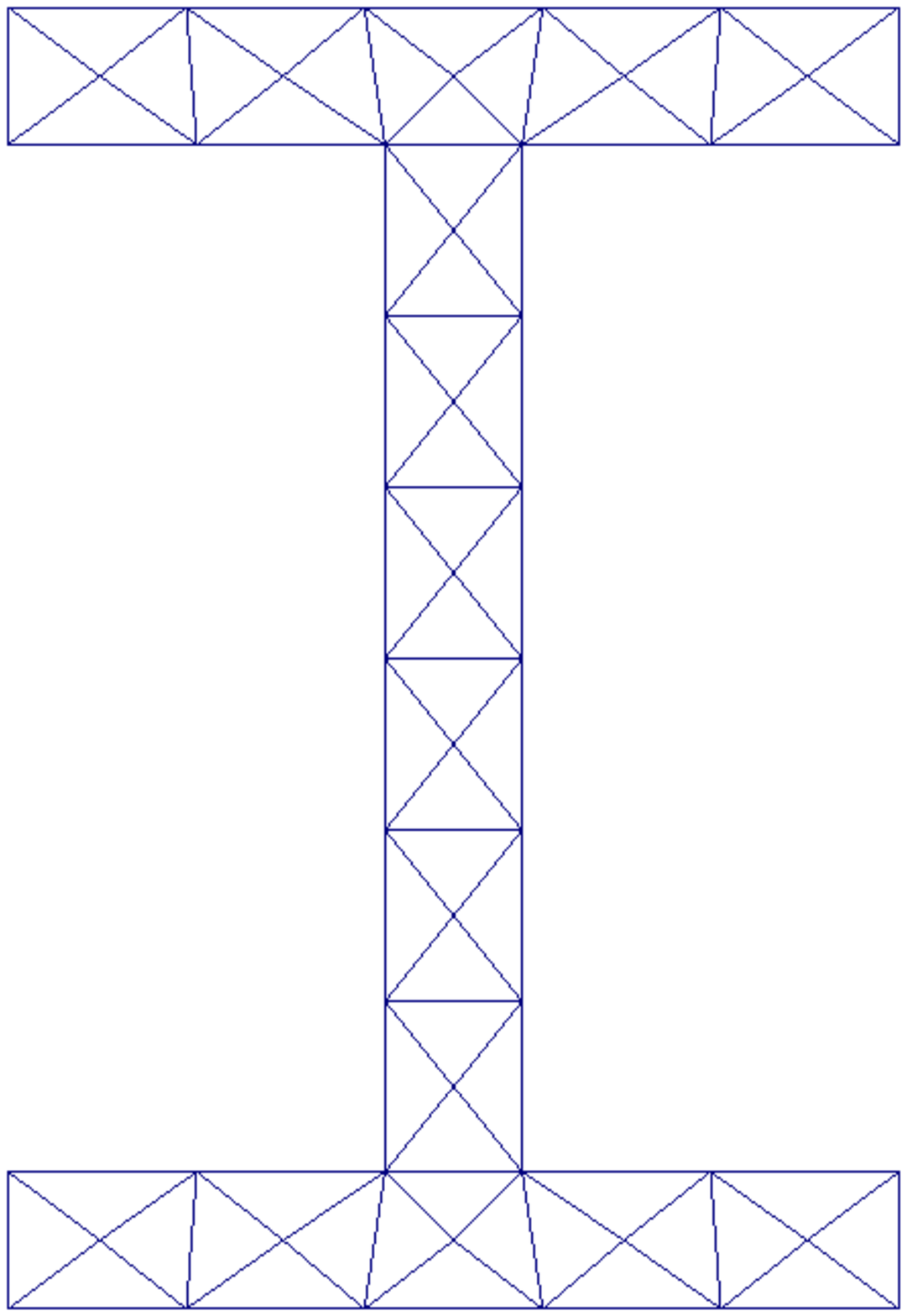} \label{fig_port_ibeam}} \quad
\end{minipage}
\caption{{\footnotesize Component mesh of i-beam (a) and i-beam with crack (b) and mesh of shared port (c).}}
\end{figure}

\begin{figure}[t]
\begin{minipage}{0.47\textwidth}
\includegraphics[scale = 0.35]{error_ibeams_diri.eps}
\caption{Relative error $\| \mathsf{u}(\Mu) - \mathsf{u}^{n+6}(\Mu)\|_{\mathcal{E}(\Mu)}/\|\mathsf{u}(\Mu)\|_{\mathcal{E}(\Mu)}$ for I-Beam with and without crack in $\Omega_{2}(\Mu)$; Dirichlet boundary conditions on outer ports of $\Omega_{1}(\Mu)$ and $\Omega_{2}(\Mu)$ are $(0,0,0)^{T}$ and $(1,1,1)^{T}$, respectively.\label{fig rel error 3d}}
\end{minipage}
\qquad
\begin{minipage}{0.47\textwidth}
\includegraphics[scale = 0.35]{eigenvalues_deviation_ibeams.eps}
\caption{Eigenvalues $\sqrt{\lambda_{n+1}(\Mu)}$ for the I-Beam with and without crack in $\Omega_{2}(\Mu)$ and deviation $E(S(R^{n}_{data,\mathrm{ker}}(\Mu)),R^m)$ during spectral greedy algorithm. {\color{white} Dirichlet boundary conditions on outer ports of $\Omega_{1}$ and $\Omega_{2}$ are $0$ and $1$}\label{fig eigenvalues 3d}}
\end{minipage}
\end{figure}

Thus, we consider next an I-Beam and a cracked I-Beam, whose corresponding component meshes are depicted in Fig.~\ref{fig_ibeam} and Fig.~\ref{fig_cracked_ibeam} and the joint port meshes can be seen in Fig.~\ref{fig_port_ibeam}. We generate the reduced port space $R^{n}_{data,\mathrm{ker}}(\Mu)$ by connecting two un-defective I-Beams. Then, we prescribe homogeneous Dirichlet boundary conditions at the outer port of $\Omega_{1}$ and $g=(1,1,1)^{T}$ at the outer port of $\Omega_{2}(\Mu)$ and assess the relative error between the FE solution $\mathsf{u}(\Mu)$ and the port reduced static condensation approximation $\mathsf{u}^{n+6}(\Mu)$ based on that port space in Fig.~\ref{fig rel error 3d}; here $+6$ accounts for the six rigid body motions included in $R^{n}_{data,\mathrm{ker}}(\Mu)$. We observe a stagnation of the relative error if we connect an I-Beam with a cracked I-Beam and use the spectral modes generated by connecting two un-defective I-Beams. However, again, we stress that we obtain a relative error of less than $10^{-3}$, which is very satisfactory. We also highlight the extremely fast convergence of the reduced static condensation approximation for the I-Beam and thus the convincing approximation capacities of the optimal port spaces for this test case.

\subsection{Spectral greedy algorithm for geometry changes}\label{subsect:num_greedy}

If we perform the spectral greedy algorithm to generate a joint port space both for the defective and un-defective I-Beam we obtain a port space of size $23$ for a tolerance of $2\cdot 10^{-6}$. Taking into account that for this tolerance the eigenspaces for the transfer eigenvalue problems for each geometry have a dimension of $16$ and $15$ including the six rigid body modes, we observe that at least for this tight tolerance neither of the two eigenspaces is well suited to approximate the other. Note, that this is consistent with our observation of the stagnation of the relative error in Fig.~\ref{fig rel error 3d}. However, based on the latter we would expect that for higher tolerances the size of the port space generated by the spectral greedy is significantly smaller than the size of space obtained by uniting the two eigenspaces obtained by the transfer eigenvalue problem. This can indeed be observed in Fig.~\ref{fig eigenvalues 3d}. If we prescribe for instance a tolerance of $10^{-2}$ the dimension of the port space obtained by the spectral greedy would be $14$ while the dimension of the union of the two eigenspaces is $17$. Increasing the tolerance further rises the gain we obtain by employing the spectral greedy rather than uniting the two eigenspaces as can be observed in Fig.~\ref{fig eigenvalues 3d}, where we compare $2\sqrt{\lambda_{n+1}(\Mu)}$ and the scaled deviation $2 E(S(R^{n}_{data,\mathrm{ker}}(\Mu)),R^m)$. Note that the factor $2$ comes from our chosen division of the tolerance in the spectral greedy, namely $\varepsilon/2$.  

\begin{figure}[t]
\begin{minipage}{0.47\textwidth}
\includegraphics[scale = 0.35]{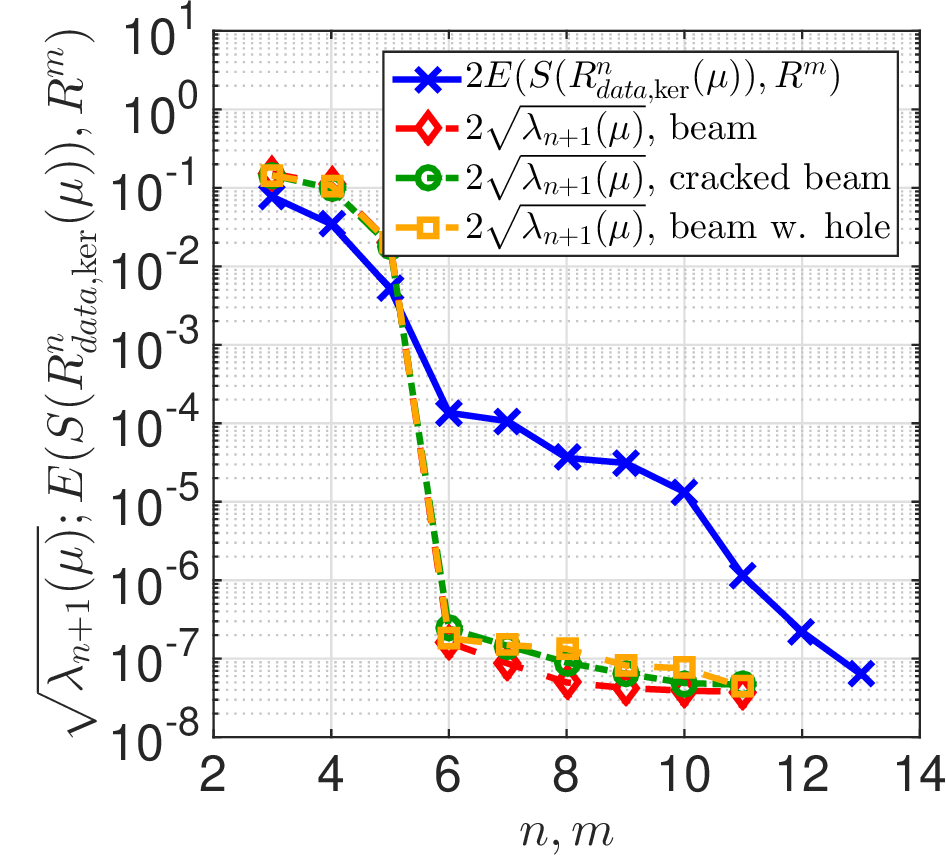}
\caption{Eigenvalues $\sqrt{\lambda_{n+1}(\Mu)}$ for the un-defective beam, the cracked beam, and the beam with a hole, and the deviation $E(S(R^{n}_{data,\mathrm{ker}}(\Mu)),R^m)$ during the spectral greedy algorithm.\label{fig dev evalues 2d}}
\end{minipage}
\qquad
\begin{minipage}{0.47\textwidth}
\includegraphics[scale = 0.35]{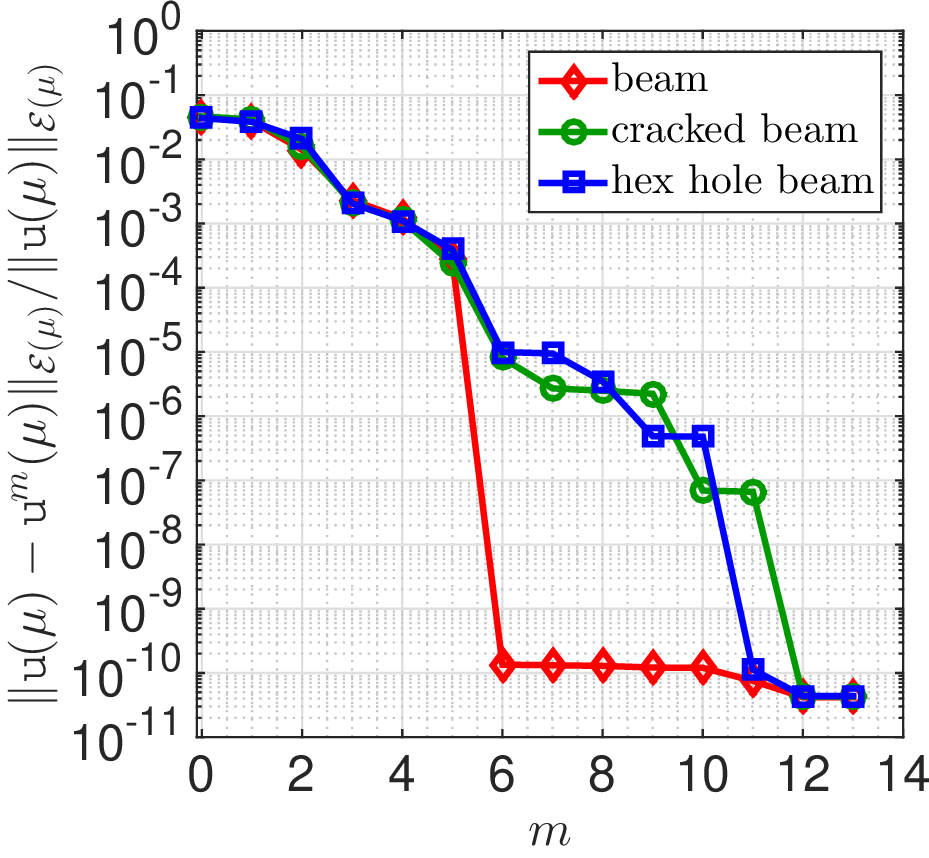}
\caption{Average relative error $\|\mathsf{u}(\Mu) - \mathsf{u}^m(\Mu)\|_{\mathcal{E}(\Mu)}/\|\mathsf{u}(\Mu)\|_{\mathcal{E}(\Mu)}$ for the un-defective beam, the cracked beam, and the beam with a hole. {\color{white}deviation $E(S(R^{n}_{data,\mathrm{ker}}(\Mu)),R^m)$ during the}\label{fig rand rel error 2d}}
\end{minipage}
\end{figure}

Similar results are obtained in two space dimensions. We connect the 2d beam as introduced in the beginning of subsection \ref{subsect:num reuse} subsequently with the 2d beam, the cracked beam depicted in Fig.~\ref{fig_no_cut_out}, and the beam with a hole (see Fig.~\ref{fig_slit_mesh}). For a tolerance of $2\cdot 10^{-7}$ the spectral greedy yields a port space of dimension $13$. As the three eigenspaces have the sizes $6$ (un-defective beam and beam with hole) and $7$, including the three rigid body modes, we observe that in this case the dimension of the port space generated by the spectral greedy equals the dimension of the union of the three eigenspaces. However, we emphasize that for larger tolerances as $10^{-3}$ or $10^{-2}$, which are of actual interest in engineering applications, we observe, again, that the spectral greedy is able to produce a very small port space which is able to yield accurate approximations for geometries which are rather different (see Fig.~\ref{fig dev evalues 2d}). 

Finally, we analyze the convergence behavior of the relative error $\| \mathsf{u}(\Mu) - \mathsf{u}^m(\Mu)\|_{\mathcal{E}(\Mu)}/\|\mathsf{u}(\Mu)\|_{\mathcal{E}(\Mu)}$ if we connect either an un-defective beam, a cracked beam, or a beam with a hole with an un-defective beam and consider random Dirichlet boundary conditions as above. Again, we observe that already for very few modes, in this case $5$, we obtain for all geometries a relative error below $10^{-3}$. However, if we insist on accuracies of $10^{-7}$ or below, we need at least for the defective components nearly all modes provided by the spectal greedy. The very good convergence behavior of the beam can be explained by the fact that after the initialization with the three rigid body modes the spectral greedy selects three (un-defective) beam modes, such that the (un-defective) beam eigenspace is contained in the spectral greedy port spaces already for $m=6$. Analyzing the convergence behavior for the defective components in detail, we observe that the modes selected from the cracked beam-beam combination reduce the error for the beam with hole-beam combination only very slightly and vice versa (see Fig.~\ref{fig rand rel error 2d}), because as the $7$th mode the cracked beam-beam, as the $8$th and $9$th the beam with hole-beam, as the $10$th the cracked beam-beam, as the $11$th the beam with hole-beam, and finally as the $12$th
 and $13$th mode the cracked beam-beam combination has been selected during the spectral greedy. 
 
\section{Conclusions}\label{sec:8}

Having the application in structural health monitoring in mind we have proposed quasi-optimal port spaces for parametrized PDEs on different geometries. To that end, we employed the optimal port spaces generated by a transfer eigenvalue problem as introduced in \cite{SmePat16} and slightly generalized the there suggested spectral greedy algorithm to geometry changes. In the numerical experiments we showed that for tolerances of actual interest the spectral greedy algorithm is able to construct a small port space that already yields an accurate approximation. Moreover, we demonstrated that using the optimal port space generated on a component pair of two un-defective beam yields a very satisfactory relative approximation error on a component pair of an un-defective beam and a beam with a defect. Inspite of the very significant change of geometry the port space can thus be reused. We therefore expect that if one constructs port spaces for a library of defects, and then detects a new defect, very often reusing the constructed port space will result in a small relative approximation error. \\

\textbf{Acknowledgements:}
I would like to thank Prof. Dr. Anthony Patera for many fruitful discussions and comments on the content of this paper. This work was supported by OSD/AFOSR/MURI Grant FA9550-09-1-0613 and ONR Grant N00014-11-1-0713.\\

\bibliographystyle{siam}

\end{document}